\documentclass[12pt, reqno]{amsproc}
\usepackage{amsmath,amssymb,mathrsfs,hyperref,color}

\DeclareMathAlphabet{\mathpzc}{OT1}{pzc}{m}{it}
\usepackage[nameinlink]{cleveref}

\usepackage{mathtools}
\usepackage[x11names]{xcolor}

\usepackage{paralist}

\usepackage{subfigure}
\usepackage{float}
\usepackage{palatino}
\usepackage{tikz}
\usepackage{cases}
\usetikzlibrary{intersections}
\usepackage{xcolor}
\usepackage[latin1]{inputenc}
\usepackage{bbm}
\usepackage{tabularx}
\usepackage{booktabs}
\usepackage[labelfont=bf,format=plain,justification=raggedright,singlelinecheck=false]{caption}
\usetikzlibrary{shadings}
\usetikzlibrary[intersections]
\usetikzlibrary{arrows,shapes}

\makeatletter
\def\setliststart#1{\setcounter{\@listctr}{#1}%
  \addtocounter{\@listctr}{-1}}
\makeatother

\makeatletter
\@addtoreset{figure}{section}
\makeatother

\usepackage{calc}

\newtheorem{remarks}{\textbf{Remark}}[section]

\newtheorem{lemma}{\textbf{Lemma}}[section]
\newtheorem{theorem}{\textbf{Theorem}}[section]
\newtheorem{corollary}{\textbf{Corollary}}[section]

\newtheorem{assumption}{\textbf{Assumption}}[section]
\numberwithin{equation}{section}

\newcommand{\R}{\mathbb{R}}

\newcommand{\PP}{\mathscr{P}}
\newcommand{\W}{\mathscr{ W}}

\newcommand{\C}{\mathcal{C}}

\DeclareMathOperator*{\ddiv}{\textrm{div}}
\DeclareMathOperator*{\ppx}{\partial_x}
\DeclareMathOperator*{\ppv}{\partial_v}
\DeclareMathOperator*{\ppy}{\partial_y}
\DeclareMathOperator*{\ppp}{\partial_p}
\DeclareMathOperator*{\ppq}{\partial_q}
\DeclareMathOperator*{\ppt}{\partial_t}

\DeclareMathOperator*{\eps}{\varepsilon}

\makeatletter
\def\moverlay{\mathpalette\mov@rlay}
\def\mov@rlay#1#2{\leavevmode\vtop{%
   \baselineskip\z@skip \lineskiplimit-\maxdimen
   \ialign{\hfil$\m@th#1##$\hfil\cr#2\crcr}}}
\newcommand{\charfusion}[3][\mathord]{
    #1{\ifx#1\mathop\vphantom{#2}\fi
        \mathpalette\mov@rlay{#2\cr#3}
      }
    \ifx#1\mathop\expandafter\displaylimits\fi}
\makeatother

\title[Rate of convergence]{Rate of convergence for first-order singular perturbation problems: Hamilton-Jacobi-Isaacs equations and mean field games of acceleration}
\author{Piermarco Cannarsa \and Cristian Mendico}
\address{Dipartimento di matematica, Universit\'a degli studi di Roma Tor Vergata -- Via della Ricerca Scientifica 1, 00133 Roma}
\email{cannarsa@mat.uniroma2.it}
\address{Dipartimento di matematica, Universit\'a degli studi di Roma Tor Vergata -- Via della Ricerca Scientifica 1, 00133 Roma}
\email{mendico@mat.uniroma2.it}

\date{\today}
\thanks{{\it Acknowledgment:} We warmly thanks the anonymous referees for the careful reading of the manuscript and the valuable comments. }

\begin{document}

\begin{abstract}
This work focuses on the rate of convergence for singular perturbation problems for first-order Hamilton-Jacobi equations. We use the nonlinear adjoint method to analyze how the Hamiltonian's regularizing effect on the initial data influences the convergence rate.  As an application we derive the rate of convergence for singularly perturbed two-players zero-sum deterministic differential games (i.e., leading to Hamilton-Jacobi-Isaacs equations) and, subsequently, in case of singularly perturbed mean field games of acceleration. Namely, we show that in both the models the rate of convergence is $\varepsilon$.

\vspace{1cm}
\begin{center}
{\it In memory of Maurizio Falcone}
\end{center}
\end{abstract}
\maketitle

\tableofcontents

\section{Introduction}

Although singular perturbation problems for first-order Hamilton-Jacobi equations have been extensively investigated over the decades, to the author' knowledge, there are no available results on the rate of convergence for such problems in the case of a non-compact state space and in the case of a fully nonlinear state equation of the underlying control problem. More specifically, we are interested in Hamilton-Jacobi-Bellman equations that arise from the dynamic programming approach to optimization, differential games, and mean field game problems with rapidly oscillating dynamics. Thus, in this work, we will consider controlled systems of the form:
\begin{equation*}
\begin{cases}
\dot x(t) = f(x(t), y(t), \alpha(t), \beta(t)), & t \in [0,T] 
\\
\dot y(t) = \frac{1}{\eps}g(x(t), y(t), \alpha(t), \beta(t)), & t \in [0,T]
\\
x(0) = x_0, \quad y(0) = y_0, & (x_0, y_0) \in \R^{d_1} \times \R^{d_2} 
\end{cases}
\end{equation*}
that lead to the following Hamilton-Jacobi-Bellman equation:
\begin{equation*}
\begin{cases}
\ppt u^{\eps}(t, x, y)  + H\left(x, y, \ppx u^{\eps}(t, x, y), \frac{1}{\eps} \ppy u^{\eps}(t, x, y)\right) = 0, & (t, x, y) \in [0,T] \times \R^{d_1} \times \R^{d_2}
\\
u^{\eps}(0, x, y) = u_{0}(x, y), & (x, y) \in \R^{d_1} \times \R^{d_2}.
\end{cases}
\end{equation*}
As $\eps \downarrow 0$, the solution of a singular perturbation problem generally leads to the elimination of the fast state variable and, consequently, to a reduction of the system's dimension. It is known (see, for instance, reference \cite{bib:ABE}) that the limit Hamilton-Jacobi equation takes the form:
\begin{equation*}
\begin{cases}
\ppt u^0(t, x) + \overline{H}(x, \ppx u^{0}(t, x)) = 0, & (t, x) \in [0,T] \times \R^{d_1}
\\
u^{0}(0, x) = \displaystyle{\min_{y \in \R^{d_2}}} u_0(x, y), & x \in \R^{d_1}
\end{cases}
\end{equation*}
where $\overline{H}$ denotes the homogenized Hamiltonian that will be introduced precisely later in \Cref{knonw_results}. Hence, as $\eps \downarrow 0$, the information on the eliminated fast variable is preserved by the Hamiltonian function.  More precisely, we show that 
\begin{equation*}
|u^{\eps}(t, x, y) - u^{\eta}(t, x, y)| \leq C_r\frac{1}{t}(C_H + t^{1/\gamma'}) (\eps - \eta)\quad \forall\; (t, x, y) \in (0,T] \times \R^{d_1} \times B_{r}.
\end{equation*}
	where $\gamma'$ is the conjugate of $\gamma$ that express the growth of the Hamiltonian w.r.t. the momentum variables. 

The theory of viscosity solutions, as presented in the seminal paper \cite{bib:LPV}, is the natural tool for this kind of problem. However, in this work, we approach viscosity solutions of first-order equations via the well-known vanishing viscosity approximation, which allows us to provide sup-norm and gradient estimates of solutions using the nonlinear adjoint method (see, for instance, \cite{Goffi_Cirant} and \cite{Evans_03} and references therein). Thus, in addition to the interest in the rate of convergence for singular perturbation problems, we believe that such a PDE method can have its own interest and provide new information on the model. In particular, we also emphasize that the assumptions on the Hamiltonian are weaker compared to those considered in the literature (e.g. \cite{Camilli_11}, \cite{Marchi_23, Marchi_11}), which are needed to solve the problem using variational techniques. Furthermore, our method addresses the lack of compactness (or, the periodicity of the coefficients), unlike the assumptions made in \cite{Italo_2001}. We also emphasize that our proposed method implicitly provides the rate of convergence for the corresponding second-order problem driven by a Laplacian operator. 

Among the many specific models that have been studied we mention, for instance, \cite{Daria} where the convergence for a class of singular perturbation problems in the whole space with the dirt of the fast variables is of Ornstein-Uhlenbeck type is addressed and, then, in \cite{Paola} such a result was extended to problems driven by a subelliptic operator. Moreover, in \cite{Cha}, \cite{Marco} ergodic problems for viscous Hamilton-Jacobi equations with superlinear growth and inward pointing drifts have been investigated and in \cite{Loreti} the connection with the long time behavior of solutions of the Cauchy problem for semi linear parabolic equations with the Ornstein-Uhlenbeck operator in the whole space is studied. Sill, there is a large amount of mathematical and engineering literature on singular perturbation problems for deterministic and stochastic control. General references for these applications include, for instance, \cite{Malley_74}, \cite{Petar_86}, \cite{Alain88}, and \cite{Yin98}. However, in addition to the engineering applications mentioned so far, in this work, we apply the results obtained for Hamilton-Jacobi-Bellman equations to singularly perturbed two-player zero-sum deterministic differential games and to mean field games with control of acceleration. Although such results have their own interest, considering the recent development of numerical methods for mean field games (especially for the Lagrangian approach, e.g., \cite{bib:Francisco1, Francisco2}), we believe that the estimation of the rate of convergence could be interesting for numerical simulations. Indeed, this rate of convergence can be used to study the convergence rates of both numerical schemes and their corresponding discrete counterparts. 

We would like to highlight that MFGs were introduced in \cite{bib:LL1, bib:LL2, bib:LL3} and \cite{bib:HCM1, bib:HCM2} to describe the behavior of Nash equilibria in problems with infinitely many rational agents. For more details on this, we refer to \cite{bib:NC} and the references therein. Since these pioneering works, MFG theory has grown rapidly. We refer, for instance, to the survey papers and monographs \cite{bib:DEV, bib:BFY, bib:CD1}. The classical MFG system describes models in which the typical payoff is represented by a deterministic calculus of variation problem. However, MFG systems with control of acceleration, first introduced in \cite{bib:CM, bib:YA}, describe models in which agents control their acceleration and the cost functional depends on higher-order derivatives of admissible trajectories. Such problems naturally arise in the study of agent-based models that describe the collective behavior of various animal populations (e.g. \cite{bib:AC, bib:TBL}) or crowd dynamics (e.g. \cite{bib:CPT, bib:EBA}). The study of the singular perturbation problem provided in this paper has numerous applications, such as in a MFG system with Cucker-Smale type dynamics (\cite{Bardi_2021}) to describe the behavior of a flock in which the control is increasingly weaker.

\medskip

\subsection*{Organization of the paper} 
 \Cref{knonw_results} presents the general model we are interested in and the known results on the convergence as the singular parameter goes to zero. In \Cref{Rate}, we present the main results of this paper, namely, we estimate the rate of convergence in the case of homogeneous Hamiltonians, and then we generalize the results to a Hamilton-Jacobi-Bellman equation arising from an optimal control problem with fully nonlinear dynamics. Finally, \Cref{Applications} focuses on the application of the previous results to differential games (see \Cref{Differential}) and to mean field games with control of acceleration (see \Cref{MFG}). \Cref{Appendix} provides sup-norm estimates and gradient estimates  that are needed in the previous parts for first-order Hamilton-Jacobi-Bellman equations.

\section{Singular perturbation problems}

In this section, we begin by introducing the general setting of the problem. We then recall the main well-known results on convergence and provide a full description of the limit equation. Finally, we present the main achievements of this work concerning the rate of convergence.

\subsection{Known convergence results}\label{knonw_results}

The following results on convergence for singular perturbation problems can be found, for instance, in \cite{bib:ABE, Bardi_ARMA}.

Let us consider the controlled dynamical system:
\begin{equation*}
\begin{cases}
\dot x(t) = f(x(t), y(t), \alpha(t)), & t \in [0,T] 
\\
\dot y(t) = \frac{1}{\eps}g(x(t), y(t), \alpha(t)), & t \in [0,T]
\\
x(0) = x_0, \quad y(0) = y_0, & (x_0, y_0) \in \R^{d_1} \times \R^{d_2} 
\end{cases}
\end{equation*}
where $\alpha: [0,T] \to \R^m$ is a measurable control function. We define the value function $u^{\eps}: [0,T] \times \R^{d_1} \times \R^{d_2} \to \R$ as
\begin{equation*}
u^{\eps}(t, x, y) = \inf_{\substack{\alpha \in L^1(0,T) \\ x(t) =x,\; y(t) =y}} \left\{ \int_{0}^{t} L(x(s), y(s), \alpha(s))\;ds + h(x, y)  \right\}
\end{equation*}
The Hamilton-Jacobi equation satisfied by $u^{\eps}$ in the viscosity sense is
\begin{equation*}
\begin{cases}
\ppt u^{\eps}(t, x, y) \\ \qquad\quad+ \displaystyle{\sup_{\alpha \in \R^m}} \left\{-\langle f(x, y, \alpha), \ppx u^{\eps}(t, x, y) \rangle - \frac{1}{\eps}\langle g(x, y, \alpha), \ppy u^{\eps}(t, x, y) - L(x, y, \alpha)\right\} = 0
\\
u(0, x, y) = h(x, y).
\end{cases}
\end{equation*}
We also define the Hamiltonian 
\[\
H(x, y, p, q) = \sup_{\alpha \in \R^m} \left\{\langle f(x, y, \alpha), p \rangle + \frac{1}{\eps}\langle g(x, y, \alpha),q \rangle - L(x, y, \alpha)\right\}.
\]
It is known that as $\eps \downarrow 0$, the limit function $\overline u$ does not depend on the fast variable $y$, and it solves a Hamilton-Jacobi equation governed by the homogenized Hamiltonian, denoted by $\overline H$. This Hamiltonian is the value of an ergodic control problem for the fast subsystem with forced slow variable $x$ and $\eps = 1$.

Thus, we can proceed as follows. Given $(\overline x, \overline p)$ consider the discounted Hamilton-Jacobi equation 
\begin{equation*}
\delta w_{\delta}(y; \overline x, \overline p) + H(\overline x, y, \overline p, \ppy w_{\delta}(y; \overline x, \overline p)) = 0, \quad y \in \R^{d_2}.
\end{equation*}
So, using classical ergodic theory, we can define
\begin{equation*}
\overline H( \overline x, \overline p) = -\lim_{\delta \to 0} \delta w_{\delta}(y; \overline x, \overline p), \quad \forall\ (\overline x, \overline p) \in \R^{2d_1}. 
\end{equation*}
For the ergodic theorem for non-compact manifolds (in this case, the Euclidean space $\R^{d_2}$), we refer to \cite[Theorem 1.1]{bib:FM}, \cite{Barles}, \cite{Ishii, Ishii2}. Therefore, as $\eps \downarrow 0$, the value function $u^{\eps}$ locally and uniformly converges to the solution $\overline u$ of the equation:
\begin{equation*}
\ppt \overline u(t, x) + \overline H(x, \ppx \overline u(t, x)) = 0, \quad (t, x) \in [0,T] \times \R^{d_1}. 
\end{equation*}

\subsection{Rate of convergence}\label{Rate}

\subsubsection{Homogenous case}

We begin our study of the rate of convergence for the singular perturbation problem by considering the case of a homogeneous Hamiltonian associated with the singular dynamical system. This first case, beyond its intrinsic interest, will prove useful later when we study the application to mean field games of acceleration.

We consider the Cauchy problem 
\begin{equation}\label{Cauchy1}
\begin{cases}
\ppt u^{\eps}(t, x, y) + \frac{1}{2\eps} |\ppy u^{\eps}(t, x, y)|^{2} + H(x, \ppx u^{\eps}(t, x, y)) = 0, & (t, x, y) \in [0,T] \times \R^{d_1} \times \R^{d_2}
\\
u^{\eps}(0, x, y) = u_{0}(x,  y), & x \in \R^{d_1}
\end{cases}
\end{equation}
where we assume 
\begin{itemize}
\item[($i$)] $H: \mathbb{R}^{d_1} \times \R^{d_1} \to \mathbb{R}$ is a Tonelli Hamiltonian, i.e., $H \in C^2(\mathbb{R}^{d_1} \times \R^{d_1})$ and the map $p \mapsto H(x, p)$ is convex and coercive;
\item[($ii$)] $u_0 : \mathbb{R}^{d_1} \times \R^{d_2} \to \mathbb{R}$ is a bounded Lipschitz continuous function.
\end{itemize} 
 From the results presented in \Cref{knonw_results}, we deduce that the viscosity solution $u^{\epsilon}: [0,T] \times \mathbb{R}^{d_1} \times \mathbb{R}^{d_2} \to \mathbb{R}$ locally uniformly converges to a continuous function $\overline{u}: [0,T] \times \mathbb{R}^{d_1} \to \mathbb{R}$ that solves the following Cauchy problem 
\begin{equation}\label{limit1} \begin{cases} 
\ppt \overline{u} (t, x) + H(x, \ppx \overline u(t, x)) = 0, & (t, x) \in [0,T] \times \R^{d_1} 
\\
\overline u(0, x, y) = \displaystyle{\min_{y \in \R^{d_2}}} u_{0}(x, y), & x \in \R^{d_1} 
\end{cases} \end{equation} 
in the viscosity sense. Indeed, the homogenized Hamiltonian for a mechanical system is the function identically equal to zero and thus we only have the contribution of the Hamiltonian $H(x, p)$ in \eqref{Cauchy1}.

Next, we prove the main result concerning system \eqref{Cauchy1}: we first determine the rate of convergence of the gradient of the solution $u^{\eps}$ w.r.t. the fast variable and, then, we provide the stability estimate w.r.t. the parameter $\eps$.

\begin{theorem}\label{Rate1}
Let $u^{\eps}$ be a solution of \eqref{Cauchy1}. Then, the following holds. 
\begin{itemize}
\item[($i$)]There exists a constant $C \geq 0$ such that 
 \begin{equation*}
  \| \ppy u^{\eps}(\tau, \cdot, \cdot)\|_{L^{\infty}(\R^{d_1} \times \R^{d_2})} \leq  \sqrt{\eps}\; C\frac{1}{\tau^{1/2}}, \quad \forall\; \tau \in (0,T].
 \end{equation*}
\item[($ii$)] For any $r \geq 0$ there exists a constant $C_r$ for any $\eps \geq \eta > 0$ we have 
\begin{equation}\label{stability1}
|u^{\eps}(t, x, y) - u^{\eta}(t, x, y)| \leq C_r \frac{1}{t^{1/2}}(\sqrt{\eps} - \sqrt{\eta})
\end{equation}
for any $(t, x, y) \in (0,T] \times \R^{d_1} \times B_{r}$.
\end{itemize}
\end{theorem}
\proof

Setting $w^{\eps}(t, x, y) = u^{\eps}(t, x, y/\sqrt{\eps})$, it solves 
\begin{equation*}
\partial_t w^{\eps}(t, x, y) + \frac{1}{2} |\ppy w^{\eps}(t, x, y)|^{2} + H(x, \ppx w^{\eps}(t, x, y)) = 0
\end{equation*}
in the viscosity sense.  Next, in order to regularize the solution, we proceed by using a vanishing viscosity approximation of $w^{\eps}$. Let $w^{\eps, \sigma}$ be a classical solution to the problem 
\begin{equation*}
\begin{cases}
\ppt w^{\eps, \sigma}(t, x, y) - \sigma \Delta w^{\eps, \sigma}(t, x, y) + \frac{1}{2} |\ppy w^{\eps, \sigma}(t, x, y)|^{2} + H(x, \ppx w^{\eps, \sigma}(t, x, y)) = 0
\\
w^{\eps}(0, x, y) = u_{0}(x, y/\sqrt{\eps}).
\end{cases}
\end{equation*}
Differentiating the equation w.r.t. $y_i$ (for $i=1, \dots, d_2)$, we deduce that $z(t, x, y) = \ppy w^{\eps, \sigma} (t, x, y)$ solves 
\begin{multline}\label{zeta}
\ppt z(t, x, y) - \sigma \Delta z(t, x, y)  + \ppy w^{\eps, \sigma}(t, x, y)\cdot \ppy z(t, x, y)
\\
+ \ppp H(x, \ppx w^{\eps, \sigma}(t, x, y)) \cdot \ppx z(t, z, y) = 0
\end{multline}
where, for simplicity of notation, we drop the subindex $i$.  Finally, the function $\widehat{z}(t, x, y) = tz(t, x, y)$ solves the equation 
\begin{multline}\label{zeta1}
\ppt \widehat{z}(t, x, y) - \sigma \Delta \widehat{z}(t, x, y)  + \ppy w^{\eps, \sigma}(t, x, y)\cdot \ppy \widehat{z}(t, x, y)
\\
+ \ppp H(x, \ppx w^{\eps, \sigma}(t, x, y)) \cdot \ppx \widehat{z}(t, z, y) =  z(t, x, y)
\end{multline}
with $\widehat{z}(0, x, y) = 0.$
Next, we show that $\widehat{z}$ is uniformly bounded w.r.t. $\sigma$ and $\eps$. Indeed, let $\rho$ be the classical solution to the Fokker-Plank equation 
\begin{equation}\label{dual1}
\begin{cases}
-\ppt \rho - \sigma \Delta \rho - \ddiv_y(\ppy w^{\eps, \sigma}(t, x, y) \rho)  - \ddiv_x(\ppp H(x, \ppx w^{\eps, \sigma}(t, x, y)) \rho) = 0
\\
\rho(\tau, x, y) = \rho_\tau(x, y)
\end{cases}
\end{equation}
where $\rho_\tau \in C^{\infty}_c(\R^{d_1} \times \R^{d_2})$ and $\|\rho_\tau\|_{L^1(\R^{d_1} \times \R^{d_2})} =1$.  Multiplying \eqref{zeta} by $\rho$ and \eqref{dual1} by $\widehat{z}$, integrating over $\R^{d_1} \times \R^{d_2}$ we get
\begin{equation*}
\int_{\R^{d_1} \times \R^{d_2}} \widehat{z}(\tau, x, y)\rho_\tau(x, y)\;dxdy \leq \int_{0}^{\tau}\int_{\R^{d_1} \times \R^{d_2}} |z(t, x, y)| \rho(t, x, y)\;dtdxdy.
\end{equation*}
Hence, by H\"older inequality and \Cref{gamma_bound} we have
\begin{multline*}
\int_{\R^{d_1} \times \R^{d_2}} \widehat{z}(\tau, x, y)\rho_\tau(x, y)\;dxdy \\ \leq \left(\int_{0}^{\tau}\int_{\R^{d_1} \times \R^{d_2}} |z(t, x, y)|^{2} \rho(t, x, y)\right)^{\frac{1}{2}} \left(\int_{0}^{\tau}\int_{\R^{d_1} \times \R^{d_2}} \rho(t, x, y)\;dtdxdy \right)^{\frac{1}{2}} \leq C \tau^{\frac{1}{2}}
\end{multline*}
which implies 
 \begin{equation*}
 \| z(\tau, \cdot, \cdot)\|_{L^{\infty}(\R^{d_1}\times \R^{d_2})} \leq C\frac{1}{\tau^{1/2}}, \quad \forall\; \tau \in (0, T].
 \end{equation*}
 Hence, passing to the limit as $\sigma \downarrow 0$ and recalling that $u^{\eps}(t, x, y) = w^{\eps}(t, x, \sqrt{\eps}y)$ we deduce 
 \begin{equation*}
  \| \ppy u^{\eps}(\tau, \cdot, \cdot)\|_{L^{\infty}(\R^{d_1} \times \R^{d_2})} \leq  \sqrt{\eps}\; C\frac{1}{\tau^{1/2}}, \quad \forall\; \tau \in (0,T].
 \end{equation*}
 Moreover, for any $\eps \geq \eta  > 0$ we have that
 \begin{multline*}
 |u^{\eps}(t, x, y) - u^{\eta}(t, x, y)| =|w^{\eps}(t, x, \sqrt{\eps}y) - w^{\eta}(t, x, \sqrt{\eta}y)| 
 \\
 \leq C|y||\ppy w^{\eps}(t, x, \sqrt{\eps}y)| (\sqrt{\eps} - \sqrt{\eta}) \leq C \frac{1}{\tau^{1/2}} |y|(\sqrt{\eps} - \sqrt{\eta})
 \end{multline*}
 which completes the proof. \qed
 
 \medskip
As an immediate consequence of the stability estimate \eqref{stability1}, we obtain the following rate of convergence of the fast solution $u^{\eps}$ to the homogenized $\overline u$.
 
 \begin{corollary}
 Let $u^{\eps}$ and $\overline u$ be a solution of \eqref{Cauchy1} and \eqref{limit1}, respectively. Then, for any $r \geq 0$ there exists a constant $C_r \geq 0$ such that
 \begin{equation*}
\|u^{\eps}(t, \cdot, \cdot) - \overline u(t, \cdot)\|_{L^{\infty}(\R^{d_1} \times B_{R})} \leq C_r\frac{1}{t^{1/2}} \sqrt{\eps}, \quad \forall\; t \in (0,T]. 
 \end{equation*}
 \end{corollary}

\begin{remarks}\em
Let us consider the problem
\begin{equation*}
\begin{cases}
\ppt u^{\eps}(t, x, y) + \frac{1}{\gamma\eps} |\ppy u^{\eps}(t, x, y)|^{\gamma} + H(x, \ppx u^{\eps}(t, x, y)) = 0, & (t, x, y) \in [0,T] \times \R^{d_1} \times \R^{d_2}
\\
u^{\eps}(0, x, y) = u_{0}(x, y), & x \in \R^{d_1}
\end{cases}
\end{equation*}
where $\gamma > 1$. Following the same reasoning as in the proof of \Cref{Rate1} we get
\begin{equation*}
|u^{\eps}(t, x, y) - u^{\eta}(t, x, y)| \leq C\frac{1}{t^{1/\gamma'}} (\varepsilon^{\frac{1}{\gamma'}} - \eta^{\frac{1}{\gamma'}})\quad \forall\; (t, x, y) \in (0,T] \times \R^{d_1} \times B_{R},\; \forall\; R > 0
\end{equation*}
	where $\gamma'$ denotes the conjugate exponent to $\gamma$. 
Consequently, 
 \begin{equation*}
\|u^{\eps}(t, \cdot, \cdot) - \overline u(t, \cdot)\|_{L^{\infty}(\R^{d_1} \times B_{R})} \leq C_R\frac{1}{t^{1/\gamma'}} \varepsilon^{\frac{1}{\gamma'}}, \quad \forall\; t \in (0,T]
 \end{equation*}
 where $\overline u$ is still the solution to the problem \eqref{limit1}. 
\end{remarks}

\begin{remarks}\em
In this setting, the initial datum of the equation depends on the fast variable. Consequently, it is natural to seek a rate of convergence with a regularizing effect. This means the rate depends on time and explodes as time approaches zero. Similar effects have been observed in \cite{Alessandro}, for example.
\end{remarks}

\subsubsection{Fully nonlinear case}\label{fully_nonlinear}

We now proceed to generalize \Cref{Rate1} to the fully nonlinear problem:
\begin{equation}\label{Cauchy2}
\begin{cases}
\ppt u^{\eps}(t, x, y)  + H\left(x, y, \ppx u^{\eps}(t, x, y), \frac{1}{\eps} \ppy u^{\eps}(t, x, y)\right) = 0, & (t, x, y) \in [0,T] \times \R^{d_1} \times \R^{d_2}
\\
u^{\eps}(0, x, y) = u_{0}(x, y), & x \in \R^{d_1}.
\end{cases}
\end{equation}

\begin{assumption}\label{H-assum}\em
As for \eqref{Cauchy1}, we assume the following. 
\begin{itemize}
\item[($i$)] $u_0: \mathbb{R}^{d_1} \times \R^{d_2} \to \mathbb{R}$ is a bounded Lipschitz continuous function.
\item[($ii$)] $H \in C^2(\R^{2d_1} \times \R^{2d_2})$ such that 
\begin{equation*}
H(x, y, p, q) \geq 0, \quad \partial^{2}_{pp}H(x, y, p, q)  > 0, \quad \partial_{qq}^{2} H(x, y, p, q) > 0, \quad \forall\; x, y \in \R^{d_1} \times \R^{d_2}
\end{equation*}
and there exist a constant $C_{H} > 0$, a constant $\gamma \geq 1$ such that for any $(x, p, y, q) \in \R^{2d_1} \times \R^{2d_2}$ 
\begin{align*}
|\ppy H(x, y, p, q)| \leq\ & C_H(1+|p|^{\gamma}+|q|^{\gamma}) \tag{H1}
\\
\langle \ppp H(x, y, p, q), p \rangle + \langle \ppq H(x, y, p, q) , q \rangle - H(x, y, p, q) \geq\ & C_H^{-1}(|p|^{\gamma} + |q|^{\gamma}) -C_H. \tag{H2}
\end{align*}
\end{itemize}
\end{assumption}

Denoting by $\overline H: \mathbb{R}^{2d_1} \to \mathbb{R}$ the homogenized Hamiltonian associated with $H$, we have that $u^\epsilon$ locally and uniformly converges to a continuous function $\overline u: [0,T] \times \mathbb{R}^d \to \mathbb{R}$ that solves:
\begin{equation}\label{limit2}
\begin{cases}
\ppt\overline u(t, x) + \overline H(x, \ppx \overline u(t, x)) = 0, & (t, x) \in [0,T] \times \R^{d_1}
\\
\overline u(0, x) =  \displaystyle{\min_{y \in \R^{d_2}}} u_{0}(x, y), & x \in \R^{d_1}
\end{cases}
\end{equation}
in the viscosity sense.

%

\medskip
\begin{theorem}\label{thm2}
Let $u^{\eps}$ be a solution of \eqref{Cauchy2}. Then, the following holds. 
\begin{itemize}
\item[($i$)] There exists a constant $C \geq 0$ such that 
\[
|\ppy u^{\eps}(t, x, y)| \leq C\frac{1}{t}(TC_H + t^{1/\gamma'}) \eps
\]
for all $(t, x, y) \in (0,T] \times \R^{d_1} \times \R^{d_2}$, where $\gamma'$ is the conjugate exponent to $\gamma$.
\item[($ii$)] For any $r \geq 0$ there exists a constant $C_r \geq 0$ such that for any $\eps \geq \eta > 0$ we have 
\begin{equation*}
|u^{\eps}(t, x, y) - u^{\eta}(t, x, y)| \leq C_r\frac{1}{t}(TC_H + t^{1/\gamma'}) (\eps - \eta)\quad \forall\; (t, x, y) \in (0,T] \times \R^{d_1} \times B_{r}
\end{equation*}
	where $\gamma'$ is the conjugate exponent to $\gamma$. 
\end{itemize}
\end{theorem}
\proof

 Define $w^{\eps}(t, x, y) = u^{\eps}(t, x, y/\eps)$ for any $(t, x, y) \in [0,T] \times \R^{d_1} \times \R^{d_2}$. Then, $w^{\eps}$ solves
\begin{equation*}
\begin{cases}
\ppt w^{\eps}(t, x, y)  + H\left(x, y/\eps, \ppx w^{\eps}(t, x, y), \ppy w^{\eps}(t, x,  y)\right) = 0, & (t, x, y) \in [0,T] \times \R^{d_1} \times \R^{d_2}
\\
w^{\eps}(0, x, y) = u_{0}(x,  y/\eps), & (x, y) \in \R^{d_1} \times \R^{d_2}
\end{cases}
\end{equation*}
in the viscosity sense. Then, proceeding as in \Cref{Rate1}, for any $\sigma > 0$ we consider the classical solution $w^{\eps, \sigma}$ to 
\begin{equation}\label{w}
\begin{cases}
\ppt w^{\eps, \sigma}(t, x, y) - \sigma \Delta w^{\eps, \sigma}(t, x, y) + H\left(x, y/\eps, \ppx w^{\eps, \sigma}(t, x, y), \ppy w^{\eps, \sigma}(t, x, y)\right) = 0, 
\\
w^{\eps}(0, x, y) = u_{0}(x, y/\eps).
\end{cases}
\end{equation}
So, differentiating the equation w.r.t. $y_i$ we deduce that the function $z^{\eps, \sigma}(t, x, y) := \ppy w^{\eps, \sigma}(t, x, y)$ is a classical solution to 
\begin{multline}\label{y_der}
\ppt z^{\eps, \sigma}(t, x, y) - \sigma \Delta z^{\eps, \sigma}(t, x, y) + \frac{1}{\eps}\ppy H\left(x, y/\eps, \ppx w^{\eps, \sigma}(t, x, y), \ppy w^{\eps, \sigma}(t, x,  y)\right) 
\\
+  \ppp H\left(x, y/\eps, \ppx w^{\eps, \sigma}(t, x,  y), \ppy w^{\eps, \sigma}(t, x, y)\right)  \cdot \ppx z^{\eps, \sigma}(t, x, y) 
\\
+  \ppq H\left(x, y/\eps, \ppx w^{\eps, \sigma}(t, x,  y), \ppy w^{\eps, \sigma}(t, x,  y)\right)  \cdot \ppy z^{\eps, \sigma}(t, x, y) = 0
\end{multline}
where, for simplicity of notation, we drop the subindex $i$. Finally, setting $\widehat{z}(t, x, y) = t z^{\eps, \sigma} (t, x, y)$ we have that $\widehat{z}$ satisfies
\begin{multline}\label{y_der1}
\ppt \widehat{z}(t, x, y) - \sigma \Delta \widehat{z}(t, x, y) + \frac{t}{\eps} \ppy H\left(x, y, \ppx w^{\eps, \sigma}(t, x, y), \ppy w^{\eps, \sigma}(t, x, y)\right) 
\\
+  \ppp H\left(x, y/\eps, \ppx w^{\eps, \sigma}(t, x,  y), \ppy w^{\eps, \sigma}(t, x,  y)\right)  \cdot \ppx \widehat{z}(t, x, y) 
\\
+  \ppq H\left(x, y/\eps, \ppx w^{\eps, \sigma}(t, x,  y), \ppy w^{\eps, \sigma}(t, x,  y)\right)  \cdot \ppy \widehat{z}(t, x, y) = z^{\eps, \sigma}(t, x, y).
\end{multline}
We consider, now, the solution $\rho$ to the dual equation 
\begin{equation}\label{rho_zeta}
\begin{cases}
-\ppt \rho(t, x, y) - \sigma\Delta \rho(t, x, y)   -  \ddiv_x(\ppp H\left(x, y/\eps, \ppx w^{\eps, \sigma}(t, x,  y),  \ppy w^{\eps, \sigma}(t, x, y)\right) \rho(t, x, y) ) \\ \qquad\qquad\qquad\qquad\qquad\qquad-  \ddiv_y(\ppq H\left(x, y/\eps, \ppx w^{\eps, \sigma}(t, x, y),  \ppy w^{\eps, \sigma}(t, x,  y)\right) \rho(t, x, y) ) = 0, 
\\
\rho(\tau, x, y) = \rho_\tau(x, y), 
\end{cases}
\end{equation}
with $\tau \in (0,T)$ fixed, $\rho_\tau \in C^{\infty}_c(\R^{d_1} \times \R^{d_2})$ and $\|\rho_\tau\|_{L^1(\R^{d_1} \times \R^{d_2})} =1$. Hence, testing \eqref{rho_zeta} w.r.t. $\widehat{z}$ and \eqref{y_der1} w.r.t. $\rho$ we obtain 
\begin{multline*}
\int_{\R^{d_1} \times \R^{d_2}} \widehat{z}(\tau, x, y)\rho_{\tau}(x, y)\;dxdy  
\\
\leq \iint_{(0,\tau) \times \R^{d_1} \times \R^{d_2}}\frac{t}{\eps} \ppy H\left(x, y, \ppx w^{\eps, \sigma}(t, x,  y), \ppy w^{\eps, \sigma}(t, x,  y)\right) \rho(t, x, y)\;dtdxdy
\\
+ \iint_{(0,\tau) \times \R^{d_1} \times \R^{d_2}} |z^{\eps, \sigma}(t, x, y)|\rho(t, x, y)\;dtdxdy.
\end{multline*}
Moreover, by {\bf {\bf (H1)}} and arguing as in \Cref{gamma_bound} applied to \eqref{w} we have
\begin{multline*}
\int_{\R^{d_1} \times \R^{d_2}} \widehat{z}(\tau, x, y)\rho_{\tau}(x, y)\;dxdy 
\\
\leq \iint_{(0,\tau) \times \R^{d_1} \times \R^{d_2}} T\varepsilon^{d_2 - 1}C_H\big(1+ |\ppx w^{\eps, \sigma}(t, x, \eps y)|^{\gamma} + |\ppy w^{\eps, \sigma}(t, x, \eps y)|^{\gamma} \big)\rho(t, x, \eps y)\;dtdxdy
\\
+  \left(\iint_{(0,\tau) \times \R^{d_1} \times \R^{d_2}} |z^{\eps, \sigma}(t, x, y)|^{\gamma}\rho(t, x, y)\;dtdxdy\right)^{\frac{1}{\gamma}} \left(\iint_{(0,\tau) \times \R^{d_1} \times \R^{d_2}}\rho(t, x, y)\;dtdxdy\right)^{\frac{1}{\gamma'}}
\end{multline*}
where we also used H\"older inequality and a change of variables by rescaling of $\eps$. We claim that there exists $C \geq 0$ such that 
\begin{equation}\label{Claim}
\iint_{(0,\tau) \times \R^{d_1} \times \R^{d_2}} \big(|\ppx w^{\eps, \sigma}(t, x, \eps y)|^{\gamma} + |\ppy w^{\eps, \sigma}(t, x, \eps y)|^{\gamma} \big)\rho(t, x, \eps y)\;dtdxdy \leq C. 
\end{equation}
If so, this yields
\begin{equation*}
\|z^{\eps, \sigma}(t, \cdot, \cdot)\|_{L^{\infty}(\R^{d_1} \times \R^{d_2})} \leq C\frac{1}{t}(TC_H\varepsilon^{d_2 - 1} + t^{1/\gamma'})
\end{equation*}
and, thus, as $\sigma \downarrow 0$ 
\begin{equation*}
\| \ppy u^{\eps}(t, \cdot, \cdot)\|_{L^{\infty}(\R^{d_1} \times \R^{d_2})} \leq C\frac{1}{t}(TC_H\varepsilon^{d_2 - 1}  + t^{1/\gamma'}) \eps \leq C\frac{1}{t}(TC_H + t^{1/\gamma'}) \eps , \quad \forall\; t \in (0,T].
\end{equation*}
In conclusion, 
 \begin{multline*}
 |u^{\eps}(t, x, y) - u^{\eta}(t, x, y)| =|w^{\eps}(t, x, \eps y) - w^{\eta}(t, x, \eta y)| 
 \\
 \leq C|y||\ppy w^{\eps}(t, x, \eps y)| (\eps - \eta) \leq C\frac{1}{t}(TC_H + t^{1/\gamma'})  |y|(\eps -\eta)
 \end{multline*}
 	for any $t \in (0,T]$.
	
	Thus, it suffices to prove \eqref{Claim}. To do so, we show that the Hamilton-Jacobi equation satisfied by $w^{\eps, \sigma}(t, x, \eps y)$ and the Fokker-Planck equation satisfied by $\rho(t, x, \eps y)$ are in duality in the sense of \eqref{lem:dualityarg} in \Cref{Appendix}. If so, \eqref{Claim} follows by the same argument in \Cref{gamma_bound}. Setting $q(t, x, v) = w^{\eps, \sigma}(t, x, \eps v)$, from \eqref{w} we deduce that $q$ is a classical solution of 
	\begin{multline}\label{HJ_q}
	\ppt q(t, x, y) - \sigma \Delta_x q(t, x, y) - \frac{\sigma}{\varepsilon^2}\Delta_y q(t, x, y)  \\+ H\left(x, y, \ppx q(t, x, y), \frac{1}{\eps}\ppy q(t, x, y)\right) = 0
	\end{multline}
	and setting $\mu(t, x, y) = \rho(t, x, \eps y)$ from \eqref{rho_zeta} we obtain that $\mu$ is a classical solution of 
	\begin{multline}\label{mu_q}
	-\ppt \mu(t, x, y) - \sigma\Delta_x \mu(t, x, y) - \frac{\sigma}{\varepsilon^2}\Delta_y \mu(t, x, y)  \\ -  \mbox{div}_x\left(\ppp H\left(x, y, \ppx q(t, x,  y),  \frac{1}{\eps} \ppy q(t, x, y)\right) \mu(t, x, y) \right) \\ -  \frac{1}{\eps}\mbox{div}_y\left(\ppq H\left(x, y, \ppx q(t, x, y),  \frac{1}{\eps} \ppy q(t, x,  y)\right) \mu(t, x, y) \right) = 0.
	\end{multline}
	Hence, we integrate and multiply \eqref{HJ_q} and \eqref{mu_q} by $\mu$ and $q$, respectively, and integrating by parts taking the sum we obtain 
	\begin{multline}\label{duality_eps}
	\int_{\R^{d_1} \times \R^{d_2}} q(\tau, x, y) \mu_{\tau}(x, y)\;dxdy - \int_{\R^{d_1} \times \R^{d_2}} q(s, x, y)\mu(s, x, y)\;dxdy 
\\
 = \iint_{(s, \tau) \times \R^{d_1} \times \R^{d_2}} \Big(\ppp H\left(x, y, \ppx q(t, x, y),  \frac{1}{\eps}\ppy q(t, x, y)\right)  \cdot \ppx q(t, x, y) 
 \\
 + \ppq H\left(x, y, \ppx q(t, x, y),  \frac{1}{\eps}\ppy q(t, x, y)\right)  \cdot \frac{1}{\eps}\ppy q(t, x, y) 
 \\
 - H\left(x, y, \ppx q(t, x, y), \frac{1}{\eps}\ppy q(t, x, y)\right)\mu(t, x, y)\;dtdxdy.
	\end{multline} 
	Finally, following the reasoning in \eqref{sup_norm} and \Cref{gamma_bound} starting from \eqref{duality_eps} we complete the proof of the claim. 
	 \qed

\medskip
Therefore, we can deduce the rate of convergence of the singularly perturbed system to the homogenized problem.

 \begin{corollary}\label{Rate2}
 Let $u^{\eps}$ and $\overline u$ be a solution of \eqref{Cauchy2} and \eqref{limit2}, respectively. Then, for any $r \geq 0$ there exists a constant $C_r \geq 0$ such that
 \begin{equation*}
\|u^{\eps}(t, \cdot, \cdot) - \overline u(t, \cdot)\|_{L^{\infty}(\R^{d_1} \times B_{r})} \leq C_r \frac{1}{t}(TC_H + t^{1/\gamma'}) \eps, \quad \forall\; t \in (0,T]. 
 \end{equation*}
 \end{corollary}


\section{Applications}\label{Applications}

\subsection{Hamilton-Jacobi-Isaacs equations}\label{Differential}

Systems, as in \eqref{Cauchy2}, arise, for instance, from deterministic two-player differential games. Consider the following systems of controlled (singularly perturbed) differential equations: 
\begin{equation*}
\begin{cases}
\dot x(t) = f(x(t), y(t), \alpha(t), \beta(t)), & t \in [0,T] 
\\
\dot y(t) = \frac{1}{\eps}g(x(t), y(t), \alpha(t), \beta(t)), & t \in [0,T]
\\
x(0) = x_0, \quad y(0) = y_0, & (x_0, y_0) \in \R^{d_1} \times \R^{d_2} 
\end{cases}
\end{equation*}
with $\alpha: [0,T] \to A$ and $\beta: [0,T] \to B$ are admissible controls with $A \subset \R^{m_1}$, $B \subset \R^{m_2}$ compact. We are also given a cost functional of the form: 
\begin{equation*}
J_t(x, y,\alpha, \beta)= \int_{0}^{t} L(x(s), y(s), \alpha(s), \beta(s))\;ds + u_0(x(0), y(0))
\end{equation*}
where $x(\cdot)$ and $y(\cdot)$ are admissible trajectories for $\alpha(\cdot)$ and $\beta(\cdot)$, respectively, with conditions $x(t) = x$ and $y(t) = y$. We assume that the first player wants to minimize the cost and the second player seeks to maximize it, i.e., we consider a singularly perturbed two-player zero-sum deterministic differential game.

Next, we define the admissible strategies and the value of the game following \cite{Fleming_89} and \cite{Evans}. First, let us denote by $\mathcal{A}(t)$ and $\mathcal{B}(t)$ the sets of admissible controls $\alpha$ and $\beta$, respectively, for the game. An admissible strategy $\alpha$ for the first player is a map $\alpha : \mathcal{A}(t) \to \mathcal{B}(t)$ such that $\alpha[\beta](\cdot)$ is an admissible control for the second player, and we define the admissible strategies for the second player in an obvious symmetric way. Moreover, we denote by $\Gamma(t)$ and $\Theta(t)$, respectively, the sets of admissible strategies for the first and second player. Finally, the well-posedness of the approach to differential game via viscosity solution to Hamilton-Jacobi equation needs the introduction of the so-called nonanticipative strategies. We recall that given $t_0 > 0$ a map $\alpha: \mathcal{B}(t_0) \to \mathcal{A}(t_0)$ is nonanticipative if for any time $t_1 > t_0$ and any controls $v_1$, $v_2 \in \mathcal{B}(t_0)$ such that $v_1(t) = v_2(t)$ a.e. $t \in [t_0, t_1]$ we have $\alpha[v_1](t) = \alpha[v_2](t)$ a.e. $t \in [t_0, t_1]$.

Then, the lower value function $u^{\eps}$ of the game is defined as
\begin{equation*}
u^{\eps}(t, x, y) = \inf_{\alpha \in \Gamma(t)} \sup_{\beta \in \mathcal{B}(t)} J(t, x, y, \alpha[\beta], \beta)
\end{equation*}
and, similarly, the upper value function $\widetilde u^{\eps}$ as
\begin{equation*}
\widetilde u^{\eps}(t, x, y) = \sup_{\beta \in \Theta(t)} \inf_{\alpha \in \mathcal{A}(t)} J(t, x, y, \alpha, \beta[\alpha]). 
\end{equation*}
We can now define the first-order Bellman-Isaacs equation starting from the upper Hamiltonian $H$ and lower Hamiltonian $\widetilde H$ as follows:
\begin{equation*}
H(x, y, p, q) = \min_{\beta \in B} \max_{\alpha \in A} \left\{-\langle p, f(x, y, \alpha, \beta) \rangle - \langle q, g(x, y, \alpha, \beta) \rangle - L(x, y, \alpha, \beta)  \right\}
\end{equation*}
and
\begin{equation*}
\widetilde H(x, y, p, q) =  \max_{\alpha \in A}\min_{\beta \in B} \left\{-\langle p, f(x, y, \alpha, \beta) \rangle - \langle q, g(x, y, \alpha, \beta) \rangle - L(x, y, \alpha, \beta)  \right\}.
\end{equation*}
Then, we have that $u^{\eps}$ is a viscosity solution to the Cauchy problem 
\begin{equation}\label{dgame1}
\begin{cases}
\partial_t u^{\eps}(t, x, y)  + H\left(x, y, \ppx u^{\eps}(t, x, y), \frac{1}{\eps} \ppy u^{\eps}(t, x, y)\right) = 0, & (t, x, y) \in [0,T] \times \R^{d_1} \times \R^{d_2}
\\
u^{\eps}(0, x, y) = u_{0}(x,  y), & (x, y) \in \R^{d_1} \times \R^{d_2}.
\end{cases}
\end{equation}
and $\widetilde u^{\eps}$ is a viscosity solution to 
\begin{equation}\label{dgame2}
\begin{cases}
\partial_t \widetilde u^{\eps}(t, x, y)  + \widetilde H\left(x, y, \ppx \widetilde u^{\eps}(t, x, y), \frac{1}{\eps} \ppy \widetilde u^{\eps}(t, x, y)\right) = 0, & (t, x, y) \in [0,T] \times \R^{d_1} \times \R^{d_2}
\\
\widetilde u^{\eps}(0, x, y) = u_{0}(x, y), & (x, y) \in \R^{d_1} \times \R^{d_2}.
\end{cases}
\end{equation}

Proceeding as before, we have that $u^{\eps}$ and $\widetilde u^{\eps}$ converge locally uniformly to $u$ and $\widetilde u$, respectively. Moreover, denoting $\overline H$ and $\widetilde H_{hom}$ as the homogenized Hamiltonians associated with $H$ and $\widetilde H$ respectively, we can conclude that $u$ and $\widetilde u$ are viscosity solutions of
\begin{equation}\label{lgame1}
\begin{cases}
\partial_t u(t, x) + \overline H(x, \ppx u(t, x)) = 0, & (t, x) \in [0,T] \times \R^{d_1}
\\
u(0, x) = \displaystyle{\min_{y \in \R^{d_2}}} u_0(x,  y), & x \in \R^{d_1}
\end{cases}
\end{equation}
and
\begin{equation}\label{lgame2}
\begin{cases}
\partial_t \widetilde u(t, x) + \widetilde H_{hom}(x, \ppx  \widetilde u(t, x)) = 0, & (t, x) \in [0,T] \times \R^{d_1}
\\
\widetilde u(0, x) = \displaystyle{\min_{y \in \R^{d_2}}} u_0(x,  y), & x \in \R^{d_1}.
\end{cases}
\end{equation}

Hence, from the same reasoning in \Cref{Rate2} we deduce the following result. 

\begin{theorem}
Let $u^{\eps}$ and $\widetilde u^{\eps}$ be solutions of \eqref{dgame1} and \eqref{dgame2}, respectively. Let $u$ and $\widetilde u$ be solutions of \eqref{lgame1} and \eqref{lgame2}, respectively. Then, for any $R \geq 0$ there exists a constant $C_R \geq 0$ such that 
 \begin{equation*}
\|u^{\eps}(t, \cdot, \cdot) -  u(t, \cdot)\|_{L^{\infty}(\R^{d_1} \times B_{R})} \leq  C_R \frac{1}{t}(C_H + t^{1/\gamma'}) \eps, \quad \forall\; t \in (0,T]
 \end{equation*}
 and
  \begin{equation*}
\|\widetilde  u^{\eps}(t, \cdot, \cdot) -  \widetilde  u(t, \cdot)\|_{L^{\infty}(\R^{d_1} \times B_{R})} \leq  C_R \frac{1}{t}(C_H + t^{1/\gamma'}) \eps, \quad \forall\; t \in (0,T].
 \end{equation*} 
\end{theorem}

\subsection{Mean field games of acceleration}\label{MFG}

We conclude this work applying the above methods and reasoning to the singular perturbation problem in case of mean field games of acceleration, i.e.,
\begin{align}\label{acc_MFG}
 \begin{cases}
 	-\ppt u^{\eps} +\frac{1}{2\eps}|\ppv u^{\eps}|^{2} - \langle \ppx u^{\eps}, v \rangle - L_{0}(x, v, \mu^{\eps}_{t})= 0, & (t,x,v) \in [0,T] \times \R^{2d}
 	\\
 	\ppt\mu^{\eps}_{t} - \langle \ppx\mu^{\eps}_{t},v \rangle - \frac{1}{\eps}\ddiv_{v}\left(\mu^{\eps}_{t}\ppv u^{\eps} \right)=0, &  (t,x,v) \in [0,T] \times \R^{2d}
\\
\mu^{\eps}_{0}=\mu_{0}, \quad u^{\eps}(T,x,v)=g(x,m^{\eps}_{T}), & (x,v) \in \R^{2d}
 \end{cases}	
 \end{align}

 \begin{assumption}\em
 We assume $L_0: \R^d \times \PP(\R^{2d}) \to \R$ to be bounded and uniformly continuous w.r.t. the measure variable, i.e., there exist a constant $C \geq 0$ and a modulus $\omega: [0, \infty) \to [0, \infty)$ such that 
 \begin{align*}
 \sup_{\substack{(x, v) \in \R^{2d} \\ \mu \in \PP(\R^{2d})}} |L_0(x, v, \mu)| \leq\ & C \tag{L1}
 \\
 |L_0(x, v,  \mu_1) - L_0(x, v, \mu_2)| \leq\ & \omega(\W_1(\mu_1, \mu_2)),  \;\; \forall\; \mu_1, \mu_2 \in \PP(\R^{2d})\tag{L2}.
 \end{align*} 
 Moreover, we assume the function $g: \R^d \times \PP(\R^d) \to \R$ be continuous and bounded w.r.t. all variables. 
 \end{assumption}
 
We recall that the space $\PP_1(\R^{2d})$ can be equipped with the Kantorovich-Rubinstein norm defined as
\begin{equation*}
\| \mu\|_0 = \sup\left\{\int_{\R^{2d}} \varphi(x, y)\;\mu(dxdy) : \varphi \in \textrm{BL}(\R^{2d})  \right\}
\end{equation*}
where $\textrm{BL}(\R^{2d})$ denotes the set of bounded Lipschitz continuous functions over $\R^{2d}$ with Lipschitz seminorm bounded by $1$. The Kantorovich-Rubinstein distance induced by the norm $\| \cdot \|_{0}$ is denoted by $\W_1$.

 It is known from \cite{Mendico} that the solution $(u^{\eps}, \mu^{\eps})$ to \eqref{acc_MFG}, where $u^{\eps}$ is a viscosity solution of the Hamilton-Jacobi equation and $\mu^{\eps}$ is a solution of the Fokker-Planck equation in the sense of distributions, converges to the pair $(u^{0}, \mu^{0})$ that solves the MFG of control system
 \begin{align}\label{MFG_control}
	\begin{cases}
		(i)\; -\ppt u^{0}(t,x) + H_{0}(x,\ppx u^{0}(t,x), \mu^{0}_{t})=0, & \quad (t,x) \in [0,T] \times \R^{d}
		\\
		(ii)\; \ppt  m^{0}_{t} - \ddiv\big( m^{0}_{t} D_p H_{0}(x,\ppx u^{0}(t,x), \mu^{0}_{t}) \big)=0, & \quad (t,x) \in [0,T] \times \R^{d}
		\\
		(iii)\; \mu^{0}_{t} = (\text{Id}(\cdot), D_p H_{0}(x,\ppx u^{0}(t,x), \mu^{0}_{t})) \sharp m^{0}_{t}, &\quad t \in [0,T]
		\\
		\mu^{0}_{0}= \mu_{0},\,\, u^{0}(T,x)=g(x, \mu^{0}_{T}), & \quad x \in \R^{d},
	\end{cases}
	\end{align}
	where
	\begin{equation*}
	H_0(x, p, m) = \sup_{v \in \R^d} \{\langle p, v \rangle - L_0(x, v, m)\},
	\end{equation*}
	$u^0$ is a viscosity solution and $\mu^0$ a solution in the sense of distribution of the limit Hamilton-Jacobi equation and of the limit Fokker-Planck equation, respectively. 
	For clearness of presentation,  we consider the system analyzed in \cite{Mendico}.

\begin{theorem}
Let $(u^{\eps}, \mu^{\eps})$ and $(u^0, \mu^0)$ be a solution of \eqref{acc_MFG} and of \eqref{MFG_control}, respectively. Then, for any $R \geq 0$ there exists a constant $C_R \geq 0$ such that for any $\eps \geq \eta > 0$ we have
\begin{equation}\label{u_osc}
|u^{\eps}(t, x, v) - u^{\eta}(t, x, v)| \leq C_R (\sqrt{\eps} - \sqrt{\eta}), \quad \forall\; (t, x, v) \in (0,T] \times \R^d \times B_R
\end{equation}
and, consequently 
\begin{equation*}
\sup_{t \in (0,T]}\|u^{\eps}(t, \cdot, \cdot ) - u^{0}(t, \cdot)\|_{L^{\infty}(\R^d \times B_R)} \leq C_R \sqrt{\eps}.
\end{equation*}
\end{theorem}

\proof 

The estimate \eqref{u_osc} can be immediately obtained following adapting the reasoning in \Cref{Rate1} to the Hamilton-Jacobi equation in \eqref{acc_MFG} since $L_0$ satisfies the required assumptions. \qed

\medskip

\begin{remarks}\em
\begin{itemize}
\item[($i$)] We observe that the rate of convergence does not depends on time since the initial datum is independent of the fast variable.
\item[($ii$)] To the best of the authors' knowledge, the above result on the rate of convergence for mean field game problems is the first of its kind in case of singular perturbation problem. Similar techniques can be applied to homogenization problems, as in \cite{bib:SLL} and \cite{bib:CDM}, when the structure of the mean field game system is preserved in the limit. 

We observe that the rate of convergence of $\mu^{\eps}$ to $\mu^0$, which depends on the rate of convergence of minimizers for the singularly perturbed action functional, remains an open problem. We consider this an interesting problem in its own right and we plan to address it in future work. A similar issue has already been raised analyzing the rate of convergence of the vanishing viscosity approximation of first-order mean field game systems in \cite{Tang}, where the authors provide the rate of convergence for the value function but not for the measure solution.
\end{itemize}
\end{remarks}

\appendix

\section{}\label{Appendix}

The fundamental tool used to provide the estimates on the solution to equation 
\begin{equation}\label{eq:HJequation}
\begin{cases}
\partial_t u(t, x, y) - \sigma\Delta u(t, x, y)  + H\left(x, y, \ppx u(t, x, y),  \ppy u(t, x, y)\right) = 0, & (t, x, y) \in [0,T] \times \R^{2d}
\\
u(0, x, y) = u_{0}(x, y), & (x, y) \in \R^{2d}.
\end{cases}
\end{equation}
is the duality method. Hence, to do so, we introduce the following Fokker-Plank equation
\begin{equation}\label{eq:FKequation}
\begin{cases}
-\ppt \rho(t, x, y) - \sigma\Delta \rho(t, x, y)   - \textrm{div}_x(\ppp H\left(x, y, \ppx u(t, x, y),  \ppy u(t, x, y)\right) \rho(t, x, y) ) \\ \qquad\qquad\qquad\qquad\qquad\qquad- \ddiv_y(\ppq H\left(x, y, \ppx u(t, x, y),  \ppy u(t, x, y)\right) \rho(t, x, y) ) = 0,
\\
\rho(\tau, x, y) = \rho_\tau(x, y),
\end{cases}
\end{equation}
for $(t, x, y) \in [0,T] \times \R^{2d}$.

In the following, we assume on $H$ the set of hypothesis given in Assumption \ref{H-assum}. Moreover, concerning the well-posedness of all the integral formulations and the integration by parts we refer to \cite[Theorem B]{Ben-Artzi} and \cite{Goffi_Cirant}.


\begin{lemma}[Integral representation formula]\label{lem:dualityarg}
Let $u$ be a solution to \eqref{eq:HJequation} and let $\mu$ be a solution to \eqref{eq:FKequation}.
Then, for any $s \in [0, \tau)$ we have that
\begin{multline}\label{eq:ABduality}
\int_{\R^{2d}} u(\tau, x, y) \rho_{\tau}(x, y)\;dxdy - \int_{\R^{2d}} u(s, x, y)\rho(s, x, y)\;dxdy 
\\
 = \iint_{(s, \tau) \times \R^{2d}} \Big(\ppp H\left(x, y, \ppx u(t, x, y),  \ppy u(t, x, y)\right)  \cdot \ppx u(t, x, y) 
 \\
 + \ppq H\left(x, y, \ppx u(t, x, y),  \ppy u(t, x, y)\right)  \cdot \ppy u(t, x, y) 
 \\
 - H\big(x, y, \ppx u(t, x, y), \ppy u(t, x, y)\Big)\rho(t, x, y)\;dtdxdy.
\end{multline}
\end{lemma}
\proof
It is enough to plug $u$ as a test function in \eqref{eq:FKequation} and conversely $\rho$ as a test function in the weak formulation of the Hamilton-Jacobi equation. The integral formula \eqref{eq:ABduality} follows by subtracting the resulting expressions.  \qed


\begin{theorem}[Sup-norm estimate]\label{sup_norm}
Let $u$ be a solution to \eqref{eq:HJequation}.
Then, 
\begin{equation}\label{eq:infinitybound}
\sup_{t \in [0,T]}\| u(t, \cdot, \cdot)\|_{L^{\infty}(\R^{2d})} \leq C
\end{equation}
for a constant $C$ independent of $\sigma > 0$. 
\end{theorem}

\proof
First, we prove a bound from above for $u$. To do so, let $\tau \in [0,T]$ and consider the solution $\mu: [0, \tau] \times \R^{2d} \to \R$ to the following problem
\begin{equation*}
\begin{cases}
\partial_{t} \mu(t, x, y) - \sigma\Delta \mu(t, x, y)= 0, & (t, x) \in [0, \tau] \times \R^{2d}
\\
\mu(\tau,x, y)=\mu_{\tau}(x, y), & x \in \R^{2d}
\end{cases}
\end{equation*}
with $\mu_{\tau} \in \C^{\infty}_{c}(\R^{2d})$. By duality arguments, i.e., proceeding as in \Cref{lem:dualityarg}, we obtain 
\begin{multline}\label{eq:duality1}
\int_{\R^{d}} u(\tau, x, y) \mu_{\tau}(x, y)\; dxdy = \int_{\R^{d}} u_{0}(x, y) \mu(0, x, u)\; dxdy 
\\ - \int_{0}^{\tau} \int_{\R^{d}} H\big(x, y, \ppx u(t, x, y), \ppy u(t, x, y)\Big)\mu(s, x, y)\; dsdxdy.
\end{multline}
Since $H(x, y, p, q) \geq 0$, $\| \mu(t, x, y)\|_{1, \R^{2d}} = 1$ for any $t \in [0, \tau]$  and $\mu\geq 0$, we get
\begin{multline*}
\int_{\R^{2d}} u_{0}(x, y) \mu(0, x, y)\ dx  \\ - \iint_{(0, \tau ) \times \R^{2d}} H\big(x, y, \ppx u(t, x, y), \ppy u(t, x, y)\Big)\mu(s, x, y)\;dxsdxdy \leq \|u_{0}(\cdot, \cdot)\|_{L^{\infty}(\R^{2d})}.
\end{multline*}
Hence, 
\begin{equation*}
\int_{\R^{2d}} u(\tau,x, y) \mu_{\tau}(x, y)\ dx \leq \|u_{0}(\cdot, \cdot)\|_{L^{\infty}(\R^{2d})}
\end{equation*}
and, thus, by passing to the supremum, over $\mu_{\tau} \geq 0$ with $\| \mu(t, x)\|_{L^1(\R^d)} = 1$, one deduces that 
\begin{equation}\label{eq:upperbound}
u(\tau, x, y) \leq \|u_{0}(\cdot, \cdot)\|_{L^{\infty}(\R^{2d})}.
\end{equation}
To prove the lower bound for $u$, we first observe that combining \eqref{eq:ABduality} and \eqref{eq:upperbound} we get 
\begin{multline*}
\iint_{(0, \tau) \times \R^{2d}} \Big(\ppp H\left(x, y, \ppx u(t, x, y),  \ppy u(t, x, y)\right)  \cdot \ppx u(t, x, y) 
 \\
 + \ppq H\left(x, y, \ppx u(t, x, y),  \ppy u(t, x, y)\right)  \cdot \ppy u(t, x, y) 
 \\
 - H\big(x, y, \ppx u(t, x, y), \ppy u(t, x, y)\Big)\rho(t, x, y)\;dtdxdy \leq 2 \|u_{0}(\cdot, \cdot)\|_{L^{\infty}(\R^{2d})}.
\end{multline*}
So, again by \eqref{eq:ABduality} and the fact that $\mu_{\tau}$ is arbitrary we get
\begin{equation}\label{eq:lowerbound}
u(\tau, x, y) \geq - \|u_{0}(\cdot, \cdot)\|_{L^{\infty}(\R^{2d})}.
\end{equation}
Therefore, \eqref{eq:upperbound} and \eqref{eq:lowerbound} yield \eqref{eq:infinitybound}. \qed

\begin{corollary}
Let $u$ be a solution of the Cauchy problem
\begin{equation*}
\begin{cases}
\partial_t u(t, x, y) + H(x, y, \ppx u(t, x, y), \ppy u(t, x, y) ) = 0, & (t, x, y) \in [0,T] \times \R^{d} \times \R^d
\\
u(0, x, y) = u_0(x, y), & (x, y) \in \R^d \times \R^d.
\end{cases}
\end{equation*}
Then, there exists a constant $C \geq 0$ such that 
\[
\sup_{t \in [0,T]} \| u(t, \cdot, \cdot)\|_{L^{\infty}(\R^{2d})} \leq C. 
\]
\end{corollary}
\proof 
The proof follows by a vanishing viscosity argument and \Cref{sup_norm}. \qed

\begin{theorem}\label{gamma_bound}
Let $u$ be a solution to \eqref{eq:HJequation} and let $\mu$ be a solution to \eqref{eq:FKequation}.
Then, there exists a constant $C \geq 0$, independent of $\sigma$, such that 
\begin{equation*}
\iint_{[0,T] \times \R^{2d}} \big(|\ppx u(t, x, y)|^{\gamma} + |\ppy u(t, x, y)|^{\gamma} \big)\rho(t, x, y)\;dtdxdy \leq C. 
\end{equation*}
\end{theorem}
\proof
The assertion follows by \Cref{lem:dualityarg}, assumption {{\bf (H2)}} and \Cref{sup_norm}. \qed

\medskip
\small{
\noindent{\bf Declarations.}

{\bf Ethical approval.} Not applicable.

{\bf Competing interests.} The authors report there are no competing interests to declare.

{\bf Authors contributions.} All authors contribueted to the research, to the editing and to the review of the manuscript.

{\bf Funding.} The author was partially supported by Istituto Nazionale di Alta Matematica, INdAM-GNAMPA project 2023, CUP E53C22001930001, by the King Abdullah University of Science and Technology (KAUST) project CRG2021-4674 ``Mean-Field Games: models, theory and computational aspects" and by the MIUR Excellence Department Project awarded to the Department of Mathematics, University of Rome Tor Vergata, CUP E83C23000330006.

{\bf Availability of data and materials.} No Data associated in the manuscript. 
}


\end{document}